\documentclass[10pt,twoside]{article}
\usepackage{hyperref}
\newcommand{\tx}[1]{\ensuremath{\text{\upshape #1}}}
\newcommand{\wt}[1]{\ensuremath{\widetilde{#1}}}
\newcommand{\se}{\ensuremath{\subseteq}}
\newcommand{\fk}{\ensuremath{\mathfrak{k}}}
\newcommand{\C}{\ensuremath{\mathds{C}}}
\newcommand{\R}{\ensuremath{\mathds{R}}}
\newcommand{\N}{\ensuremath{\mathds{N}}}
\newcommand{\from}{\ensuremath{\nobreak:\nobreak}}
\renewcommand{\to}{\ensuremath{\nobreak\rightarrow\nobreak}}
\usepackage{amssymb,amsmath,amscd,amsthm,dsfont} 
\usepackage{graphicx}
\usepackage{hyperref}
\usepackage{calc}

\theoremstyle{definition}
\newtheorem{definition}{Definition}[section]
\newtheorem{remark}[definition]{Remark}
\theoremstyle{plain}
\newtheorem{lemma}[definition]{Lemma}
\newtheorem{proposition}[definition]{Proposition}
\newtheorem{theorem}[definition]{Theorem}
\newtheorem{corollary}[definition]{Corollary}

\newenvironment{prf}{\begin{proof}[\textbf{\upshape Proof.}]}{\end{proof}}

\usepackage[a4paper,hscale=0.7,vscale=0.8,hmarginratio=1:1,includehead,headheight=14pt]{geometry}
\usepackage{fancyhdr}
\fancypagestyle{plain}{\fancyhead{}\fancyfoot{}}

\newcommand{\shortTitle}{}

\begin{document}
\sloppy
\thispagestyle{plain} 
\title{\textbf{Smooth Extensions and Spaces of Smooth and Holomorphic 
       Mappings}\\[12pt] 
       \large - preprint -}
\renewcommand{\shortTitle}{Smooth Extensions and Spaces of
Smooth and Holomorphic Mappings}
\author{Christoph Wockel\\
        Fachbereich Mathematik\\
        Technische Universit\"at Darmstadt\\
\small \texttt{wockel@mathematik.tu-darmstadt.de}}
\maketitle
\thispagestyle{empty}
\begin{abstract}
\noindent In this paper we present another notion of a smooth manifold
with corners and relate it to the commonly used concept in the
literature.  Afterwards we introduce complex manifolds with corners
show that if $M$ is a compact (respectively complex) manifold with
corners and $K$ is a smooth (respectively complex) Lie group, then
$C^{\infty}(M,K)$ (respectively $C^{\infty}_{\C}(M,K)$) is a smooth
(respectively complex) Lie group.\\[\baselineskip]
\end{abstract}
\label{first}
\section{Introduction}We introduce the notion of a smooth manifold
with corners, which is an extension of the existing notion of smooth
manifolds with corners or boundary for the finite-dimensional case (cf.
\cite{lee03} or \cite[Chapter 2]{michor80}). The notation presented
here is the appropriate notion for a treatment of mapping spaces and
Whitney's extension theorem \cite{whitney34} implies that for
finite-dimensional smooth manifolds our definition coincides with the
one given in \cite{lee03}. We give an alternative proof of a similar
statement by elementary methods from real Analysis (cf. also
\cite[Theorem 22.17]{krieglmichor97} and \cite[Proposition
24.10]{krieglmichor97}).

Afterwards we introduce complex manifolds with corners and show
several properties of the spaces $C^{\infty}(M,K)$ and
$C^{\infty}_{\C}(M,K)$.  Eventually it turns out that these mapping
spaces are smooth (respectively complex) Lie groups. This is in particular
interesting since it seems to be the only way to put a complex or
even smooth structure onto spaces of holomorphic mappings since very
elementary results for the non-boundary case imply that compactly
supported holomorphic maps are constant%
. With the results of this paper a
topological treatment of spaces like $C^{\infty}_{\C}(M,K)$ for
non-compact $M$ becomes possible by inductive limit methods.
\section{Notions of Differential Calculus}\label{sect:diffclac}
In this section we present the elementary notions of differential
calculus on locally convex spaces and for not necessarily open
domains.
\begin{definition}\label{diffcalc1}
Let $E$ and $F$ be a locally convex spaces and $U\se E$ be
open. Then $f\from U\to F$ is \textit{continuously differentiable} or
\textit{$C^{1}$} if it is continuous, for each $v\in E$ the
differential quotient
\[
df (x).v:=\lim_{h\to 0}\frac{f (x+hv)-f (x)}{h}
\]
exists and if the map $df\from U\times E\to F$ is continuous.  If $n>1$ we
inductively define $f$ to be \textit{$C^{n}$} if it is $C^{n-1}$ and
$d^{n-1}f$ is $C^{1}$ and to be $C^{\infty}$ or \textit{smooth} if it is
$C^{n}$ for all $n\in\N$.  We denote the corresponding set of maps by
$C^{1}(U,E)$, $C^{n}(U,E)$ and $C^{\infty}(U,E)$.  This is the notion
of differentiability used in %
\cite{gloeckner02a} and it will be the notion throughout this paper.

If $E$ and $F$ are complex vector spaces, then $f$ is called
\textit{holomorphic} if it is $C^{1}$ and the map $df(x)\from E\to F$
is complex linear for all $x\in U$ (cf. \cite[p.1027]{milnor83}).
\end{definition}
\begin{remark}
From the above definition it is clear what the notions of a smooth
(respectively complex) Lie group is, i.e. a group which is a smooth
(respectively complex) manifold modelled on a locally convex (complex) space
such that the group operations are smooth (respectively holomorphic).
\end{remark}
\begin{remark}\label{convenientcalc}
(cf. \cite[Remark 3.2]{neeb03}) We briefly recall the basic
definitions underlying the convenient calculus from
\cite{krieglmichor97}. Again let $E$ and $F$ be locally convex spaces.
A curve $f\from \R\to E$ is called smooth if it is smooth in the sense of
Definition \ref{diffcalc1}.  Then the $c^{\infty}$-topology on $E$ is
the final topology induced from all smooth curves $f\in C^{\infty}
(\R,E)$. If $E$ is a Fr\'echet space, then the $c^{\infty}$-topology
is again a locally convex vector topology which coincides with the
original topology \cite[Theorem 4.11]{krieglmichor97}.  If $U\se E$ is
$c^{\infty}$-open then $f\from U\to F$ is said to be $C^{\infty}$
or smooth if
\[
f_{*}\left(C^{\infty} (\R,U) \right)\se C^{\infty} (\R,F),
\]
e.g. if $f$ maps smooth curves to smooth curves. The chain rule
\cite[Proposition 1.15]{gloeckner02b} implies that each smooth map in
the sense of Definition \ref{diffcalc1} is smooth in the convenient
sense. On the other hand \cite[Theorem 12.8]{krieglmichor97} implies
that on a Fr\`echet space a smooth map in the convenient sense is
smooth in the sense of Definition \ref{diffcalc1}. Hence for Fr\'echet
spaces the two notions coincide.
\end{remark}
\begin{definition}\label{diffcalc2}
Let $E$ and $F$ be a locally convex space, and let $U\se E$ be a set
with dense interior. We say that a map $f\from U\to F$ is \textit{$C^{1}$}
if it is continuous, $f_{\tx{int}}:=\left.f \right|_{\tx{int}(U)}$ is
$C^{1}$ and the map $d(f_{\tx{int}})$
extends to a continuous map on $U\times E$, which is called the
\textit{differential} $df$ of $f$. If $n>1$ we inductively define $f$
to be \textit{$C^{n}$} if if is $C^{1}$ and $df$ is $C^{n-1}$ for
$n>1$. We say that $f$ is $C^{\infty}$ or \textit{smooth} if $f$ is
$C^{n}$ for all $n\in \N_{0}$.
\end{definition}
\begin{remark}\label{smoothnessremark}
Since $\tx{int}(U\times E^{2^{n}-1})=\tx{int}(U)\times E^{2^{n}-1}$ we have for
$n=1$ that $\left(df \right)_{\tx{int}}=d\left(f_{\tx{int}} \right)$ and we
inductively obtain
$\left(d^{n}f \right)_{\tx{int}}=d^{n}\left(f_{\tx{int}} \right)$. Hence the
higher differentials $d^{n}f$ are defined to be the continuous extensions of
the differentials $d^{n}(f_{\tx{int}})$ and thus we have that a map
$f\from U\to F$ is smooth if and only if
$d^{n}(f_{\tx{int}})$
has a continuous extension $d^{n}f$ to $U\times E^{2^{n}-1}$ for all
$n\in \N$.
\end{remark}
\begin{remark}\label{chainruleforboundaries}
If $f\from U\to U'$, $g\from U'\to F$ with $f(\tx{int}(U))\se\tx{int}(U')$ are
$C^{1}$, then the chain rule for locally convex spaces
\cite[Proposition 1.15]{gloeckner02b} and $(g\circ
f)_{\tx{int}}=g_{\tx{int}}\circ f_{\tx{int}}$ imply that
$g\circ f\from U\to F$ is $C^{1}$ and its differential is given by
$d(g\circ f)(x).v=dg(f(x)).df(x,v)$.
In particular it follows that $g\circ f$ is smooth if $g$ and $f$ are so.
\end{remark}
\begin{definition}\label{manifoldwithboundary}
(cf. \cite{lee03} for the
finite-dimensional case) Let $E$ be a locally convex space,
$\lambda_{1},\dots ,\lambda_{n}$ be continuous functionals and 
$E^{+}:=\bigcap_{k=1}^{n}\lambda_{k}^{-1}(\R_{0}^{+})$. If $M$ is a
Hausdorff space, then a collection $(U_{i},\varphi_{i})_{i\in I}$ of
homeomorphisms $\varphi_{i}\from U_{i}\to \varphi (U_{i})$ called charts onto open
subsets $\varphi_{i}(U_{i})$ of $E^{+}$ is a
\textit{differential structure} on $M$ of co-dimension $n$ if
$\cup_{i\in I}U_{i}=M$ and for each pair of charts
$\varphi_{i}$ and $\varphi_{j}$ with $U_{i}\cap U_{j}\neq \emptyset$ we have
that the \textit{coordinate change}
\[
\varphi_{i} \left(U_{i}\cap U_{j}\right)\ni x \mapsto
\varphi_{j}\left(\varphi_{i}^{-1}(x)\right) \in \varphi_{j} (U_{i}\cap
U_{j})
\]
is smooth in the sense of Definition \ref{diffcalc2}. Two differential
structures are called \textit{compatible} if their union is again a
differential structure, a maximal differential structure with respect
to compatibility is called an atlas and $M$ together with an atlas is
called a \textit{smooth manifold with corners} of co-dimension $n$.
\end{definition}
\begin{remark}\label{usualnotion}
Note that the previous definition of a smooth manifold with corners coincides
for $E=\R^{n}$ with the one given in \cite{lee03} and in the case of
co-dimension $1$ and a Banach space $E$ with the definition of a
manifold with boundary in \cite{lang99}, but our notion of smoothness
differs. In both cases a map $f$ defined on a non-open subset $U\se E$
is said to be smooth if for each point $x\in U$ there exists an open
neighbourhood $V_{x}\se E$ of $x$ and a smooth map $f_{x}$ defined on
$V_{x}$ with $f=f_{x}$ on $U\cap V_{x}$. However, it will turn out
that for finite-dimensional smooth manifolds with corners the two notions
coincide.
\end{remark}
\begin{lemma}\label{invarianceofinteriorpoints}
If $M$ is a smooth manifold with corners modelled on the locally
convex space $E$ and $\varphi_{i}$ and $\varphi_{j}$ are two charts
with $U_{i}\cap U_{j}\neq \emptyset$, then
\mbox{$\varphi_{j}\circ\varphi_{i}^{-1}(\tx{int}(\varphi_{i}
(U_{i}\cap U_{j})))\se \tx{int}(\varphi_{j}(U_{i}\cap U_{j}))$}.
\end{lemma}
\begin{prf}
Denote by $\alpha\from \varphi_{i}(U_{i}\cap U_{j})\to
\varphi_{j}(U_{i}\cap U_{j})$, $x\mapsto
\varphi_{j}(\varphi_{i}^{-1}(x))$ and $\beta =\alpha^{-1}$ the
corresponding coordinate changes.  We claim that $d\alpha (x)\from E\to E$
is an isomorphism if $x\in\tx{int}(\varphi_{i} (U_{i}\cap U_{j}))$. Since
$\beta$ maps a neighbourhood $W_{x}$ of $\alpha (x)$ into
$\tx{int}(\varphi_{i} (U_{i}\cap U_{j}))$ we have $d\alpha (\beta
(x')).\big(d\beta (x').v\big)=v$ for $v\in E$ and $x'\in
\tx{int}(W_{x})$ (cf. Remark \ref{chainruleforboundaries}). Since
$x',v\mapsto d\alpha (\beta (x')).\big(d\beta (x').v\big)$ is
continuous and $\tx{int}(W_{x})$ is dense in $W_{x}$, we thus have
that $v\mapsto d\beta (\alpha (x)).v$ is a continuous inverse for
$d\alpha (x)$.

Now suppose $x\in \tx{int}(\varphi_{i} (U_{i}\cap U_{j}))$ and $\alpha
(x)\notin \tx{int}(\varphi_{j} (U_{i}\cap U_{j}))$. Then
$\lambda_{i}(\alpha (x))=0$ for some $i\in\{1,\dots ,n\}$ and thus
there exists an $v\in E$ such that $\alpha (x)+t v\in \varphi_{j}
(U_{i}\cap U_{j})$ for $t\in [0,1]$ and $\alpha (x)+t v\notin
\varphi_{j} (U_{i}\cap U_{j})$ for $t\in [-1,0)$. But then
$v\notin\tx{im}(d\alpha (x))$, contradicting the surjectivity of
$d\alpha (x)$.
\end{prf}
\begin{remark}
The preceding Lemma shows that the points of
$\tx{int}(E_{+})$
are invariant under coordinate changes and thus the \textit{interior}
$\tx{int}(M)=\bigcup_{i\in I}\varphi_{i}^{-1}(\tx{int}(E_{+}))$
is an intrinsic property of $M$. We denote by $\partial
M:=M\backslash \tx{int}(M)$ th boundary of $M$.
\end{remark}
\begin{definition}\label{smoothonmanifold}
A map $f\from M\to N$ between smooth manifolds with corners is said to be
\textit{$C^{n}$}, respectively \textit{smooth}, if
$f\left(\tx{int}(M)\right)\se \tx{int}(N)$ and the
corresponding \textit{coordinate representation}
\[
\varphi_{i} (U_{i}\cap f^{-1}(U_{j}))\ni
x\mapsto \varphi_{j}\left(f\left(\varphi_{i}^{-1}(x)\right)\right)
\in \varphi_{j}(U_{j})
\]
is $C^{n}$ respectively smooth for each pair $\varphi_{i}$ and $\varphi_{j}$
of charts on $M$ and $N$. We again
denote the corresponding spaces of mappings by $C^{n}(M,N)$ and
$C^{\infty}(M,N)$.
\end{definition}
\begin{remark}
For a map $f$ to be smooth it suffices to check that
\[
\varphi(U\cap f^{-1}(V))\ni x\mapsto
\varphi'(f(\varphi^{-1}(x)))\in \psi (V)
\]
maps $\tx{int}(\varphi (U\cap f^{-1}(V)))$ into $\tx{int}\psi (V)$ and is
smooth in the sense of Definition \ref{diffcalc2}.
for each $m\in M$ and an arbitrary pair of charts $\varphi\from U\to E^{+}$
and $\psi \from V\to F_{+}$ around $m$ and $f(m)$ due to Lemma
\ref{chainruleforboundaries} and Lemma
\ref{invarianceofinteriorpoints}.
\end{remark}
\begin{definition}
If $M$ is a smooth manifold with corners and differentiable structure
$\left(U_{i},\varphi_{i} \right)_{i\in I}$, which is modelled on the
locally convex space $E$, then the \textit{tangent space} in
$m\in M$ is defined to be $TM_{m}:=\left(E\times I_{m}\right)/ \sim$,
where $I_{m}:=\{i\in I:m\in U_{i}\}$ and $(x,i)\sim
\left(d\left(\varphi_{j}\circ\varphi_{i}^{-1}\right)
(\varphi_{i}(m)).x,j\right)$. The set $TM:=\cup_{m\in M}\{m\}\times
T_{m}M$ is called the \textit{tangent bundle} of $M$.
\end{definition}
\begin{remark}
Note that the tangent spaces $T_{m}M$ are isomorphic for all $m\in M$,
including the boundary points.
\end{remark}
\begin{proposition}
The tangent bundle $TM$ is a smooth manifold with corners and the map
$\pi\from TM\to M$, $\left(m,[x,i]\right)\mapsto m$ is smooth.
\end{proposition}
\begin{prf}
Fix a differentiable structure $(U_{i},\varphi_{i})_{i\in I}$ on $M$.
Then each
$U_{i}$ is a smooth manifold with corners with respect to the differential
structure $(U_{i},\varphi_{i})$ on $U_{i}$.
We endow each $TU_{i}$ with the topology induced from the mappings
\begin{align*}
\tx{pr}_{1} & \from TU_{i}\to M,\;\; \left(m,v\right)\mapsto m\\
\tx{pr}_{2} & \from TU_{i}\to E,\;\; \left(m,v\right)\mapsto v,
\end{align*}
and endow $TM$ with the topology making each map $TU_{i}\hookrightarrow  TM$,
$(m,v)\mapsto (x,[v,i])$ a topological embedding. Then $\varphi_{i} \circ
\tx{pr}_{1}\times \tx{pr}_{2}\from TU_{i}\to \varphi (U_{i})\times E$
defines a differential structure on $TM$ and from the very definition
it follows immediately that $\pi$ is smooth.
\end{prf}
\begin{corollary}\label{smoothmapsonmanifolds}
If $M$ and $N$ are smooth manifolds with corners, then a map $f\from M\to N$ is
$C^{1}$ if $f(\tx{int}(M))\se \tx{int}(N)$,
$f_{\tx{int}}:=\left.f \right|_{\tx{int}(M)}$ is $C^{1}$ and
$Tf_{\tx{int}}\from T(\tx{int}(M))\to T(\tx{int}(N))\se TN$ extends
continuously to $TM$. If, in addition, $f$ is $C^{n}$ for
$n\geq 2$, then the map
\[
Tf\from TM\to TN,\;\;\left(m,[x,i] \right)\mapsto
\left(f(m),[d\left(\varphi_{j}\circ f\circ\varphi_{i}^{-1}\right)
		\left(\varphi_{i} (m)\right).x,j] \right)
\]
is well-defined and $C^{n-1}$.
 \end{corollary}
\begin{definition}
If $M$ is a smooth manifold with corners, then for $n\in\N_{0}$ the
\textit{higher tangent bundles} $T^{n}M$ are the inductively defined
smooth manifolds with corners \mbox{$T^{0}M:=M$} and
\mbox{$T^{n}:=T\left(T^{n-1}M\right)$}. If $N$ is a smooth manifold
with corners and $f\from M\to N$ is $C^{n}$, then the \textit{higher
tangent maps} $T^{m}f\from T^{m}M\to T^{m}N$ are the inductively defined
maps $T^{0}f:=f$ and $T^{m}f:=T(T^{m-1}f)$ if $1<m\leq n$.
\end{definition}
\begin{corollary}
If $M$, $N$ and $O$ are smooth manifolds with corners and $f\from M\to N$ and
$g\from N\to O$ with $f(\tx{int}(M))\se \tx{int}(N)$ and $g(\tx{int}(N))\se
\tx{int}(O)$ are $C^{n}$, then $f\circ g\from M\to O$ is
$C^{n}$ and we have $T^{m} (g\circ f)=T^{m}f\circ T^{m}g$ for all
$m\leq n$.
\end{corollary}

\begin{proposition}\label{partitionsofunity}
If $M$ is a finite-dimensional paracompact smooth manifold with
corners and $(U_{i})_{i\in I}$ is an open cover of $M$, then there
exists a smooth partition of unity $(f_{i})_{i\in I}$ subordinated to
this open cover.
\end{proposition}
\begin{prf}
The construction in \cite[Theorem 2.1]{hirsch76} actually yields smooth
functions $f_{i}\from U_{i}\to \R$ also in the sense of Definition
\ref{smoothonmanifold}.
\end{prf}
\section{Extensions of Smooth Maps}\label{extensionofsmoothmaps}
This section draws on a suggestion by Helge
Gl{\"o}ckner and was inspired by \cite[Chapter IV]{broecker}.  We
relate the notions introduced in Definition \ref{diffcalc2} to the
usual notion of differentiability on a non-open subset $U\se\R^{n}$
(cf. Remark \ref{usualnotion}).
We will see that, at least under some
mild requirements, this notion coincides with the definition
given in Definition \ref{diffcalc2}.

It should be emphasised that the results of this secrtion are not new (cf.
\cite{whitney34}, \cite[Theorem 22.17]{krieglmichor97} and
 \cite[Theorem 24.10]{krieglmichor97}). The remarkable thing is that it is
proved by methods of elementary Analysis.
\begin{remark}\label{rem:topologyOnSpacesOfSmoothMappings}
If $M$ is a smooth manifold with corners, we endow $C^{\infty}(M,F)$ with the
topology making the canonical map
\[
C^{\infty}(M,F)\hookrightarrow
\prod_{n\in \N_{0}}C(T^{n}M,F)_{c},\;\;f\mapsto d^{n}f
\]
a topological embedding (cf. \cite[Definition 3.1]{gloeckner02a}), where
$d^{n}f:=pr_{2^{n}}\circ T^{n}f$ (note $T^{n}F=F^{2^{n}}$).
This is a locally convex vector topology on $C^{\infty}(M,F)$ and it
is complete whenever $C^{\infty}(\tx{int}(M),F)$ is complete. Hence it
is a Fr\'echet space if $M$ is a second countable finite-dimensional
smooth manifold with corners and $F$ is a Fr\'echet space%
. Note that this is \textit{not} immediate
if one uses the notion of smoothness on $M$ from \cite{lee03} or
\cite{lang99}.
\end{remark}
\begin{definition}\label{def:CCP}
We say that the two locally convex spaces $E$ and $E'$ satisfy the
\textit{Cartesian closedness principle} (shortly CCP) if for each pair of
open subsets $U\se E$ and $U'\se E'$ and each locally convex space $F$
we have a well defined and linear isomorphism
\[
~^{\wedge}\from C^{\infty}(U\times U',F)\to C^{\infty}(U,C^{\infty}(U',F)),\;\;
f^{\wedge}(x)(y)=f(x,y).
\]
\end{definition}
\begin{remark}\label{rem:CCP}
Since for Fr\'echet spaces, our notion of differentiability and the
one used in the convenient setting \cite{krieglmichor97} coincide,
any two Fr\'echet spaces satisfy the Cartesian closedness principle
due to \cite[Lemma 3.12]{krieglmichor97}.
\end{remark}
\begin{proposition}\label{prop:CCP}
If $E$ and $E'$ satisfy the Cartesian closedness principle, $U\se
E$, $U'\se E'$ have dense interior and $F$ is locally convex, then we
have a well-defined canonical map
\[
~^{\wedge}\from C^{\infty}(U\times U',F)\to
C^{\infty}(U,C^{\infty}(U',F)),\;\;
f^{\wedge}(x)(y)=f(x,y).
\]
If, moreover, $E'$ is finite-dimensional, then this map is a linear
isomorphism.
\end{proposition}
\begin{prf}
First we check that $f^{\wedge}$ actually is an element of
$C^{\infty}(U,C^{\infty}(U',F))$. If $x\in \tx{int}(U)$, then
$\left.f^{\wedge}(x)\right|_{\tx{int}(U')}\in
C^{\infty}(\tx{int}(U'),F)$ due to the Cartesian closedness principle
for $\tx{int}(U)$ and $\tx{int}(U')$. Since $d^{n}f$ extends
continuously to the boundary so does $d^{n}(f^{\wedge}(x))$. So
$\!\!~^{\wedge}$ defines a map from $\tx{int}(U)$ to
$C^{\infty}(U',F)$, which is continuous since $(C(X\times
Y,Z))^{\wedge}\se C(X,C(Y,Z))$. With Remark \ref{smoothnessremark}, the
smoothness of $f^{\wedge}$ follows in the same way as the continuity.

It is immediate that $\!\!~^{\wedge}$ is
linear and injective.
If $E'$ is finite-dimensional, then surjectivity follows directly from the
the Cartesian closedness principle for $E$ and $E'$ and
\mbox{$C(X\times Y,Z)\cong C(X,C(Y,Z))$}.
\end{prf}
\begin{lemma}
If $E$ is a locally convex space and $\left(f_{n}\right)_{n\in \N_{0}}$ is
a sequence in $C^{1}(\R,E)$ such that $\left(f'_{n} \right)_{n\in \N_{0}}$
converges uniformly on compact subsets to some \mbox{$\bar{f}\in C(\R,E)$},
then $\left(f_{n} \right)$ converges to some $f\in C^{1}(\R,E)$ with
$f'=\bar{f}$.
\end{lemma}
\begin{prf}
This can be proved as in the case $E=\R$ (cf. \cite[Proposition
IV.1.7]{broecker}).
\end{prf}
\begin{lemma}
Let $F$ be a Fr\'echet space. If $(v_{n})_{n\in \N_{0}}$ is an
arbitrary sequence in $F$, then there exists an $f\in
C^{\infty}(\R,F)$ such that $f^{(n)}(0)=v_{n}$ for all $n\in\N_{0}$.
\end{lemma}
\begin{prf}
(cf. \cite[Proposition IV.4.5]{broecker} for the case $F=\R$).  Let
$\zeta\in C^{\infty}(\R,\R)$ be such that $\tx{supp}(\zeta)\se[-1,1]$ and
$\zeta(x)=1$ if $-\frac{1}{2}\leq x\leq \frac{1}{2}$ and put
$\xi(x):=x\,\zeta(x)$. Then $\tx{supp}(\xi)\se[-1,1]$ and $\left.\xi
\right|_{[-\frac{1}{2},\frac{1}{2}]}=
\tx{id}_{[-\frac{1}{2},\frac{1}{2}]}$.  Since $\xi^{k}$ is compactly
supported, there exists for each $n\in \N$ an element $M_{n,k}\in\R$
such that $|\left(\xi^{k}\right)^{(n)}(x)|\leq M_{n,k}$ for all
$x\in\R$.  Now let $(p_{m})_{m\in \N}$ be a sequence of seminorms
defining the topology on $F$ with $p_{1}\leq p_{2}\leq \ldots$. We now
choose $c_{k}>1$ such that
$|p_{k}(v_{k})\,c_{k}^{-k}\left(\xi(c_{k}\,\,\cdot\,)^{k}\right)^{(n)}(x)|\leq
p_{k}(v_{k})%
c_{k}^{n-k}\,M_{n,k}<2^{-k}$
if $n<k$. Note that this is possible since there are only finitely
many inequalities for each $k$. Set
$f_{m}:=\sum_{k=0}^{m}v_{k}\left(c_{k}^{-1}\xi(c_{k}\,\,\cdot\,)\right)^{k}$.
We show that $f:=\lim_{m\to\infty}f_{m}$ has the desired properties.  If
$\varepsilon >0$ and $l\in \N$ we let $m_{\varepsilon ,l}> l$ be such
that $2^{-m_{\varepsilon ,l}}<\varepsilon$. Thus
\begin{multline*}
p_{l}(f_{m}^{(n)}-f_{m_{\varepsilon ,l}}^{(n)})=
p_{l}\big(\sum_{k=1+m_{\varepsilon ,l}}^{m} v_{k}c_{k}^{-k}
(\xi(c_{k}\;\cdot\; )^{k})^{(n)}\big)\\
\leq \sum_{k=1+m_{\varepsilon ,l}}^{m}p_{k}(v_{k})%
c^{n-k}_{k}
M_{n,k}
\leq 2^{-m_{\varepsilon,l}} < \varepsilon
\end{multline*}
for all $m>m_{\varepsilon ,l}$ and $n<l$.  It follows for $n<l$ that
$f^{(n)}_{m}$ converges uniformly to some \mbox{$f^{n}\in C^{\infty}(\R,F)$}
and the preceding lemma implies $(f^{n-1})'=f^{n}$, whence
$f^{(n)}=f^{n}$. Since $l$ was chosen arbitrarily if follows that $f$
is smooth.  Since we may also interchange differentiation and the
limit by the preceding lemma and since $c_{k}\xi(c_{k}\,\,\cdot\,)$
equals the identity on a zero neighbourhood we thus have
$f^{(n)}(0)=\left(\lim_{_{m\to \infty }}f_{m}^{(n)}\right)(0)=
\lim_{_{m\to \infty }}\left(f_{m}^{(n)}(0) \right)=v_{n}$.
\end{prf}
\begin{corollary}\label{extensionToR}
If $F$ is a Fr\'echet space then for each $f\in C^{\infty}\left([0,1],F\right)$
there exists a \mbox{$\bar{f}\in C^{\infty}(\R,F)$} with
$\left.\bar{f} \right|_{[0,1]}=f$.
\end{corollary}
\begin{prf}
(cf. \cite[Proposition 24.10]{krieglmichor97})
For $n\in\N_{0}$ set $v_{n}:=f^{(n)}(0)$ and $w_{n}:=f^{(n)}(1)$. Then
the preceding lemma yields $f_{-},f_{+}\in C^{\infty}(\R,F)$ with
$f^{(n)}_{-}(0)=v_{n}=f^{(n)}(0)$ and
\mbox{$f_{+}^{(n)}(0)=w_{n}=f^{(n)}(1)$}. Then
\[
\bar{f}(x):=
\left\{\begin{tabular}{ll}
	$f_{-}(x)$&	if $x<0$ \\
	$f(x)$& 	if $0\leq x\leq 1$ \\
	$f_{+}(x-1)$&	if $x>1$
\end{tabular} \right.
\]
defines a function on $\R$ which has continuous differentials of arbitrary
order and hence is smooth.
\end{prf}
\begin{theorem}\label{extensionTHM}
If $F$ is a Fr\'echet space and $f\in C^{\infty}([0,1]^{n},F)$, then there
exists an $\bar{f}\in C^{\infty}(\R^{n},F)$ with
$\left.\bar{f} \right|_{[0,1]^{n}}=f$.
\end{theorem}
\begin{prf}
First we note that the Cartesian closedness principle holds in this
context due to Remark \ref{rem:CCP} and
Proposition \ref{prop:CCP}.  Set $f_{0}:=f$. Using
Proposition \ref{prop:CCP} we can view $f_{0}$ as an element
\[
f_{0}\in C^{\infty}\left([0,1],C^{\infty}\left([0,1]^{n-1},E \right) \right).
\]
which we can extend to an element of
$C^{\infty}\left(\R,C^{\infty}\left([0,1]^{n-1},E \right) \right)$ 
by Corollary \ref{extensionToR} and Remark
\ref{rem:topologyOnSpacesOfSmoothMappings}.
This can again be seen as an element $f_{1}\in
C^{\infty}\left(\R\times[0,1]^{n-1},E \right)$. In the same manner we
obtain a map
\[
f_{2}\in C^{\infty}\left(\R^{2}\times[0,1]^{n-2},E \right)
\]
extending $f_{1}$ as well as $f_{0}$. Iterating this procedure for
each argument results in a map $\bar{f}:=f_{n}$ which extends each
$f_{i}$ and so it extends $f_{0}=f$.
\end{prf}
\begin{proposition}\label{extensiononmanifolds}
If $F$ is a Fr\'echet space, $M$ is a finite-dimensional smooth
manifold and $L\se M$ has dense interior and is a smooth manifold with
corners with respect to the charts obtained from the restriction of
the charts of $M$ to $L$, then there exists an open subset $U\se M$
with $L\se U$ such that for each $f\in C^{\infty}(L,F)$ there exists a
$\bar{f}\in C^{\infty}(U,F)$ with $\left.\bar{f} \right|_{L}=f$.
\end{proposition}
\begin{prf}
For each $m\in\partial L\cap L$ there exists a set $L_{m}$ which is open in
$M$
and a chart $\varphi_{m} \from L_{m}\to \R^{n}$ such that
$\varphi_{m} (L\cap L_{m})\se \R_{+}^{n}$ and
$\varphi_{m}(m)\in \partial \R_{+}^{n}$. Then there exists a cube
\[
C_{m}:=[x_{1}-\varepsilon,x_{1}+\varepsilon]\times\ldots\times
       [x_{n}-\varepsilon,x_{n}+\varepsilon]\se \varphi_{m} (L\cap L_{m}),
\]
where
\[
x_{i}=\left\{\begin{tabular}{ll}
$\varphi_{m} (m)_{i}$ & if $\varphi_{m} (m)_{i}\neq 0$\\
$\varepsilon $    & if $\varphi_{m} (m)_{i}=0$
\end{tabular}\right.
\]
(actually $C_{m}$ is contained in $\R_{+}^{n}$ and shares the $i$-th
``boundary-face'' with $\R_{+}^{n}$ if \mbox{$\varphi_{m} (m)_{i}=0$}). 
Then $C$ is
diffeomorphic to $[0,1]^{n}$. The diffeomorphism is defined by multiplication
and addition and extends to a diffeomorphism of $\R^{n}$. We now set
$U=\tx{int}(L)\cup \bigcup_{m\in \partial L\cap L}V_{m}$,
$V_{m}:=\tx{int}(\varphi_{m}^{-1}(C_{m}))$
and choose a partition of unity $g,h,(h_{m})_{m\in
\partial L\cap L}$ subordinated to the open cover $U\backslash
L,\tx{int}(L),(V_{m})_{m\in \partial L\cap L}$.

If $f\in C^{\infty}(L,M)$, then Theorem \ref{extensionTHM} yields a
smooth extension $f_{m}$ of $\left.f\circ \varphi_{m}^{-1}\right|_{C_{m}}$
and thus $\bar{f}_{m}:=\left.f_{m}\circ \varphi\right|_{V_{m}}$ is smooth
and extends $f$.  We now set
\[
\bar{f}(x):=h(x)\,f(x)+\sum_{m\in \partial L}h_{m}(x)\,f_{m}(x),
\]
where we extend $f$ and $f_{m}$ by zero if not defined. Since $h$
(respectively $h_{m}$) vanishes on a neighbourhood of each point in
$\partial L\cap L$ (respectively $\partial V_{m}$), this function is smooth
and
since $\left.f_{m} \right|_{V_{m}\cap L}=\left.f \right|_{V_{m}\cap
L}$ for all $m\in\partial L$ it also extends $f$.
\end{prf}
\begin{corollary}
If $U\se (\R^{n})^{+}$ is open, $F$ a Fr\'echet space and $f\from U\to F$ is
smooth in the sense of Definition \ref{diffcalc2}, then there exists an open
subset $\wt{U}\se \R^{n}$, with $U\se \wt{U}$, such that for each
$f\in C^{\infty}(U,F)$ there exists a $\wt{f}\in C^{\infty}(\wt{U},F)$
with $\left.\wt{f}\right|_{U}=f$.
\end{corollary}

\section{Spaces of Mappings}
In this section we prove several results on mapping spaces like
$C^{\infty}(M,K)$ or $C^{\infty}_{\C}(M,K)$.
Since many proofs carry over from the non-boundary case, we provide here
only the necessary changes to the statements and extensions to
\cite[p.366-375]{gloeckner02a}.
\begin{definition}\label{def:holomorphicMap}
If $E$ and $F$ are locally convex complex vector spaces, $U\se E$ has
dense interior, then a smooth map $f\from U\to F$ is called
\textit{holomorphic} if $f_{\tx{int}}$ is holomorphic, i.e. that
$df_{\tx{int}}(x)\from E\to F$ is complex linear
(cf. \cite[p. 1027]{milnor83}). We denote the space of all holomorphic
functions on $U$ by $C^{\infty}_{\C}(U,F)$.
\end{definition}
\begin{remark}
Note that in the above setting $df(x)$ is complex linear for all $x\in U$ due
to the continuity of the extension of $df_{\tx{int}}$.
\end{remark}
\begin{definition}
A smooth manifold with corners is called a \textit{complex manifold
with corners} if it is modelled on a complex vector space $E$ and the
coordinate changes in Definition \ref{manifoldwithboundary} are
holomorphic. A map smooth $f\from M\to N$ between complex manifolds with
corners is said to be holomorphic if and for each pair of charts on
$M$ and $N$ the corresponding coordinate representation is holomorphic
(cf. Definition \ref{smoothonmanifold}).  We denote the space of
holomorphic mappings from $M$ to $N$ by $C^{\infty}_{\C}(M,N)$.
\end{definition}
\begin{remark}
If $M$ is a complex manifold with corners and $F$ is a locally convex
complex vector space, then $C^{\infty}_{\C}(M,F)$ is a closed subspace
of $C^{\infty}(M,F)$ since the requirement on $df(x)$ being complex
linear is a closed condition as an equational requirement on $df(x)$
in the topology defined in Remark
\ref{rem:topologyOnSpacesOfSmoothMappings}.
\end{remark}
\begin{lemma}\label{teclem3}
If $M$ is a finite-dimensional smooth manifold with corners and $E$ and $F$
are locally convex spaces, then there is an isomorphism
\mbox{$C^{\infty}(M,E\times F)\cong C^{\infty}(M,E)\times C^{\infty}(M,F)$}.
\end{lemma}

\begin{prf}
The proof of \cite[Lemma 3.4]{gloeckner02a} carries over without changes.
\end{prf}

\begin{lemma}
If $M$ and $N$ are finite-dimensional smooth manifolds with corners,
$E$ is locally convex and $f\from N\to M$ is smooth, then the map
$C^{\infty} (M,E)\to C^{\infty} (N,E)$, $\gamma \mapsto \gamma \circ f$
is continuous.
\end{lemma}

\begin{prf}
The proof of \cite[Lemma 3.7]{gloeckner02a} carries over without changes.
\end{prf}

\begin{lemma}\label{teclem1}
If $M$ is a finite-dimensional smooth manifold with corners and $E$ is a
locally convex space, then the map
$C^{\infty}(M,E)\to C^{\infty}(T^{n}M,T^{n}E)$, $\gamma \mapsto T^{n}\gamma$
is continuous.
\end{lemma}

\begin{prf}
The proof of \cite[Lemma 3.8]{gloeckner02a} carries over for $n=1$, where
\cite[Lemma 3.7]{gloeckner02a} has to be substituted by Lemma \ref{teclem1}
and
\cite[Lemma 3.4]{gloeckner02a} has to be substituted by Lemma \ref{teclem3}.
The assertion follows from an easy induction.
\end{prf}

\begin{lemma}\label{teclem2}
If $X$ is a Hausdorff space, $E$ and $F$ are locally convex spaces,
$U\se E$ is open and $f\from X\times U\to F$ is continuous, then the mapping
\[
f_{\sharp}\from C(X,U)_{c}\to C(X,F)_{c},\;\;
\gamma \mapsto f\circ( \tx{id}_{X},\gamma)
\]
is continuous.
\end{lemma}

\begin{prf}
Since the topology of compact convergence and the compact-open
topology coincide on $C(X,E)$ and $C(X,F)$ \cite[Theorem
X.3.4.2]{bourbakiTop}, this is \cite[Lemma 3.9]{gloeckner02a}.
\end{prf}

\begin{lemma}\label{ppfiscont}
If $M$ is a finite-dimensional smooth manifold with corners, $E$ and $F$ are
locally convex spaces, $U\se E$ is open and $f\from M\times U\to F$ is
smooth, then the mapping
\[
f_{\sharp}\from C^{\infty}(M,U)\to C^{\infty}(M,F),\;\;\gamma\mapsto f\circ
(\tx{id}_{M},\gamma)
\]
is continuous.
\end{lemma}
\begin{prf}
For $\gamma \in C^{\infty} (M,U)$ we have
\[
T (f_{\sharp}\gamma)=T (f\circ (\tx{id}_{M},\gamma))=Tf\circ T
(\tx{id}_{M},\gamma)=Tf\circ (\tx{id}_{TM},T\gamma)= (Tf)_{\sharp} (T\gamma )
\]
and thus inductively
\begin{multline*}
T^{n} (f_{\sharp}\gamma)=T\big(T^{n-1}(f_{\sharp}\gamma)\big)=
T\big((T^{n-1}f)_{\sharp}T^{n-1}\gamma\big)\\
=T\big(T^{n-1}f\circ(\tx{id}_{T^{n-1}M},T^{n-1}\gamma)\big)=
T^{n}f\circ(\tx{id}_{T^{n}M},T^{n}\gamma)=
\big(T^{n}f\big)_{\sharp}T^{n}\gamma.
\end{multline*}
Now we can write the map $\gamma \mapsto T^{n} (f_{\sharp}\gamma)$ as the
composition of the two maps 
\mbox{$\gamma \mapsto (\tx{id}_{T^{n}M},T^{n}\gamma)$}
and
$(\tx{id}_{T^{n}M},T^{n}\gamma) \mapsto (T^{n}f)_{\sharp}T^{n}\gamma$ which
are
continuous by Lemma \ref{teclem1} and \ref{teclem2}. Hence $f_{\sharp}$
is continuous because a map from any topological space to $C^{\infty} (M,F)$
is
continuous if all compositions with $d^{n}=\tx{pr}_{2^{n}}\circ T^{n}$ are
continuous.
\end{prf}
\begin{proposition}\label{ppfissmooth}
If $M$ is a compact smooth manifold with corners, $E$ and $F$
are locally convex spaces, $U\se E$ is open and $f\from M\times U\to F$ is smooth,
then the mapping $f_{\sharp}\from C^{\infty}(M,U)\to C^{\infty}(M,F)$, $
\gamma \mapsto f\circ (\tx{id}_{M},\gamma)$
is smooth.  If, moreover, $E$ and $F$ are complex vector spaces and
$f(m)\from U\to F$ is holomorphic for all $m\in M$, then $f_{\sharp}$ is
holomorphic.
\end{proposition}
\begin{prf}
(cf. \cite[Proposition III.7]{neeb01b}) We claim that
\begin{equation}\label{dpf}
d^{n} (f_{\sharp})= (d^{n}_{2}f)_{\sharp}
\end{equation}
holds for all $n\in\mathbb{N}_{0}$, where
$d^{n}_{2}f(x,y).v:=d^{n}f(x,y).(0,v)$.
This claim immediately proves the assertion due to Lemma \ref{ppfiscont}.

To verify (\ref{dpf}) we perform an induction on $n$. The case $n=0$ is
trivial, hence assume that \eqref{dpf} holds for $n\in\mathbb{N}_{0}$ and take
\[
\gamma\in C^{\infty}(M,U)\times C^{\infty}(M,E)^{2^{n}-1}\cong
C^{\infty}(M,U\times E^{2^{n}-1})
\]
and
\[
\eta \in C^{\infty} (M,E)^{2^{n}}\cong C^{\infty}(M,E^{2^{n}}).
\]
Then $\text{im}(\gamma)\se U\times E^{2^{n}-1}$ and
$\text{im} (\eta )\se E^{2^{n}}$
are compact and there exists an $\varepsilon > 0$ such that
\[
\text{im}(\gamma)+(-\varepsilon,\varepsilon)\text{im}(\eta)\se
U\times E^{2^{n}-1}.
\]
Hence $\gamma+h\eta\in C^{\infty}(M,U\times E^{2^{n}-1})$ for all
$h\in (-\varepsilon,\varepsilon)$ and we calculate
\begin{align*}
\big(d(d^{n}f_{\sharp})(\gamma,\eta)\big)(x)&~=\lim_{h\to 0}\frac{1}{h}
\Big(\big(d^{n}f_{\sharp}(\gamma+h\eta)-d^{n}f_{\sharp}(\gamma)\big)(x)\Big)\\
&~\stackrel{i)}{=}\lim_{h\to 0}\frac{1}{h}\Big(d_{2}^{n}f
\big(x,\gamma(x)+h\eta(x)\big)-d_{2}^{n}f\big(x,\gamma(x)\big)\Big)\\
&~\stackrel{ii)}{=}\lim_{h\to 0}\int_{0}^{1}d_{2}
\bigg(\Big(d_{2}^{n}f\big(x,\gamma(x)+th\,\eta(x)\big)\Big),\eta(x)\bigg)dt\\
&~\stackrel{iii)}{=}\int_{0}^{1}\lim_{h\to 0}d_{2}
\bigg(\Big(d_{2}^{n}f\big(x,\gamma(x)+th\,\eta(x)\big)\Big),\eta(x)\bigg)dt\\
&~=d_{2}^{n+1}f\big(x,\gamma(x),\eta(x)\big)=
\big(d^{n+1}_{2}f\big)_{\sharp}(\gamma,\eta)(x),
\end{align*}
where $i)$ holds by the induction hypothesis, $ii)$ holds by the Fundamental
Theorem of Calculus \cite[Theorem 1.5]{gloeckner02b} and $iii)$ holds
due to the differentiability of parameter-dependent Integrals
(cf. \cite{gloecknerneeb}).
The derived formula
$d(f_{\sharp})=(d_{2}f)_{\sharp}$ shows that $d(f_{\sharp})$ is
complex linear.
\end{prf}
\begin{corollary}\label{pfissmooth}
If $M$ is a compact smooth manifold with corners, $E$ and $F$ are locally
convex spaces, $U\se E$ are open and $f\from U\to F$ is smooth
(respectively holomorphic), then the push forward
$f_{*}\from C^{\infty}(M,U)\to C^{\infty}(M,F)$, $\gamma \mapsto f\circ
\gamma$ is a smooth (respectively holomorphic) map.
\end{corollary}
\begin{prf}
Define $\tilde{f}\from M\times U\to F$, $(x,v)\mapsto f (x)$ and apply Proposition
\ref{ppfissmooth}.
\end{prf}
\begin{proposition}\label{localdescriptionofliegroups}
Let $G$ be a group with a smooth manifold structure on
$U\se G$ modelled on the locally convex space $E$.
Furthermore assume that there exists $V\se U$ open such that $e\in V$,
$VV\se U$, $V=V^{-1}$ and
\begin{itemize}
\item [i)]   $V\times V\to U$, $(g,h)\mapsto gh$ is smooth,
\item [ii)]  $V\to V$, $g\mapsto g^{-1}$ is smooth,
\item [iii)] for all $g\in G$ there exists an open unit neighbourhood $W\se U$
             such
             that $g^{-1}Wg\se U$ and the map $W\to U$, $h\mapsto g^{-1}hg$ is
             smooth.
\end{itemize}
Then there exists a unique smooth manifold structure on $G$ such that $V$ is
an open
submanifold of $G$ which turns $G$ into a Lie group.
\end{proposition}

\begin{prf}
The proof of \cite[Proposition III.1.9.18]{bourbakiLie} carries
over without changes.
\end{prf}
\begin{theorem}
Let $M$ be a compact smooth manifold with corners, $K$ be a Lie group and let
$\varphi \from W\to \varphi (W)\se \fk :=L(K)$ be a chart of $K$ around $e$
with $\varphi (e)=0$. Furthermore let $\varphi_{*}\from C^{\infty}
(M,W)\to C^{\infty} (M,\fk )$, $\gamma \mapsto \varphi \circ \gamma$.
\begin{itemize}
\item [a)] If $M$ and $K$ are smooth, then $\varphi_{*}$ induces a
	   smooth manifold structure on $C^{\infty}(M,K)$, turning it
	   into a smooth Lie group w.r.t. pointwise operations.
\item [b)] If $M$ is smooth and $K$ is complex, then $\varphi_{*}$ induces a
	   complex manifold structure on $C^{\infty}(M,K)$, turning it
	   into a complex Lie group w.r.t. pointwise operations.
\item [c)] If $M$ and $K$ are complex, then the restriction of $\varphi_{*}$
	   to $C^{\infty}_{\C}(M,W)$ induces a complex manifold structure on
	   $C^{\infty}_{\C}(M,K)$, turning it into a complex Lie group w.r.t.
	   pointwise operations, modelled on $C^{\infty}_{\C}(M,\fk)$.
\end{itemize}
\end{theorem}
\begin{prf}
Using Corollary \ref{pfissmooth} and Proposition \ref{ppfissmooth},
the proof of the smooth case in \cite[Section 3.2]{gloeckner02a}
carries over to yield a).  Since Proposition \ref{ppfissmooth} also
implies the holomorhpy of the group operations, b) is now immediate.
Using the same argumentation as in a) we deduce c) since $\varphi_{*}$
maps $C^{\infty}_{\C}(M,W)$ bijectively to $C^{\infty}_{\C}(M,\varphi
(W))$, which is open in $C^{\infty}_{\C}(M,\fk)$.
\end{prf}

\section*{Acknowledgements} The work on this paper was financially
supported by a doctoral scholarship from the Land Hessen. The author
would also like thank Karl-Hermann Neeb and Helge Gl\"ockner for
giving several hints to the results of this paper.

\bibliographystyle{myalpha}
\bibliography{mybib}

\vskip\baselineskip
\vskip\baselineskip
\vskip\baselineskip
\large
\noindent
Christoph Wockel\\
Fachbereich Mathematik\\
Technische Universit\"at Darmstadt\\
Schlossgartenstrasse 7\\
D-64289 Darmstadt\\
Germany\\[\baselineskip]
\normalsize
\texttt{wockel@mathematik.tu-darmstadt.de}

\label{last}
\end{document}